\documentclass[a4paper]{amsart}
\usepackage{verbatim}
\usepackage[dvips]{color}
\usepackage{amssymb}

\author[F. Duzaar]{Frank Duzaar}
\address{Frank Duzaar\\Department Mathematik, Universit\"at
Erlangen--N\"urnberg\\ Bismarckstrasse 1 1/2, 91054 Erlangen, Germany}
\email{duzaar@mi.uni-erlangen.de}

\author[G. Mingione]{Giuseppe Mingione}
\address{Dipartimento di Matematica, Universit\`a di Parma\\
Viale G.~P.~Usberti 53/a, Campus, 43100 Parma, Italy}
\email{giuseppe.mingione@unipr.it.}

\allowdisplaybreaks

\newtheorem{theorem}{Theorem}[section]

\theoremstyle{definition}

\numberwithin{equation}{section}

\newcommand{\potm}[1]{{\bf I}_{#1}^{|\mu|}}

\newcommand{\rif}[1]{(\ref{#1})}

\newcommand\ap{``}

\newcommand{\pot}[1]{{\bf I}_{#1}^{\mu}}

\def\eqn#1$$#2$${\begin{equation}\label#1#2\end{equation}}
\def\charfn_#1{{\raise1.2pt\hbox{$\chi
_{\kern-1pt\lower3pt\hbox{{$\scriptstyle#1$}}}$}}}

\def\qq1{q_*}
\def\q2{q_{**}}

\def\ep{\varepsilon}
\def\en{\mathbb N}
\def\er{\mathbb R}

\def\loc{\operatorname{loc}}

\newdimen\vintbar
\def\vint{-\kern-\vintbar\int}

\def\0{\boldsymbol 0}

\newcommand{\ratio}{\nu, L}

\newcommand{\divo}{\textnormal{div}}

 \newcommand{\mean}[1]{-\hskip-1.08em\int_{#1}}

\newcommand{\npma}{{\bf W}_{\frac{1}{p},p}^{\mu}}

\newcommand{\ww}{{\bf W}_{1,p}^{\mu}}

\newcommand{\trif}[1] {\textnormal{\rif{#1}}}

\newtoks\by
\newtoks\paper
\newtoks\book
\newtoks\jour
\newtoks\yr
\newtoks\pages
\newtoks\vol
\newtoks\publ
\def\et{ \& }
\def\name[#1, #2]{#1 #2}
\def\ota{{\hbox{\bf ???}}}
\def\cLear{\by=\ota\paper=\ota\book=\ota\jour=\ota\yr=\ota
\pages=\ota\vol=\ota\publ=\ota}
\def\endpaper{\the\by, \textit{\the\paper},
{\the\jour} \textbf{\the\vol} (\the\yr), \the\pages.\cLear}
\def\endbook{\the\by, \textit{\the\book},
\the\publ, \the\yr.\cLear}
\def\endpap{\the\by, \textit{\the\paper}, \the\jour.\cLear}
\def\endproc{\the\by, \textit{\the\paper}, \the\book, \the\publ,
\the\yr, \the\pages.\cLear}

\title[Gradient estimates in non-linear potential theory] {Gradient estimates in non-linear potential theory}
\begin{document}
Rendiconti Lincei - Matematica e Applicazioni, 20 (2009) 179-190
\maketitle
\centerline{{\em To the memory Renato Caccioppoli, Neapolitan Mathematician}}

\begin{abstract}
We present pointwise gradient bounds for solutions to $p$-Laplacean type non-homogeneous equations employing non-linear Wolff type potentials, and then prove similar bounds, via suitable caloric potentials, for solutions to parabolic equations. A method of proof entails a family of non-local Caccioppoli inequalities, together with a DeGiorgi's type fractional iteration.
\end{abstract}
\section{The classical setting and a zero order estimate}
In this note we describe some of the results and techniques developed in the papers \cite{DDM, mis3}, which give a complete non-linear analog of the classical pointwise gradient estimates valid for the Poisson equation
\eqn{poisson}
$$
-\triangle u = \mu\qquad \qquad \mbox{in}\ \er^n\,,
$$
where $\mu$ is in the most general case a Radon measure with finite total mass. Moreover, the estimates we present hold for non-linear parabolic equations. At the same time our results give a somehow unexpected but natural maximal order - and parabolic - version of a by now classical result due to Kilpel\"ainen \& Mal\'y \cite{KM} and later extended, by mean of a different approach, by Trudinger \& Wang \cite{TW}. To better frame our setting, let us recall a few basic linear results concerning the basic example \rif{poisson}
- here for simplicity considered in the whole $\er^n$ - for which, due to the use of classical representation formulas, it is possible to get pointwise bounds for solutions via the use of Riezs potentials
\eqn{ellr}
$$
I_\beta(\mu)(x):=\int_{\er^n}\frac{d \mu(y)}{|x-y|^{n-\beta}}\,,\qquad \qquad \beta \in (0,n]
$$
such as
\eqn{stima0}
$$
|u(x)|\leq  c I_2(|\mu|)(x)\,, \qquad \mbox{and}\qquad
|Du(x)|\leq  c I_1(|\mu|)(x)\,.
$$
We recall that the equivalent, localized version of the Riesz potential $I_\beta(\mu)(x)$
is given by the linear potential
\eqn{riebeta}
$$
\pot{\beta}(x_0,R):= \int_0^R \frac{\mu(B(x_0,\varrho))}{\varrho^{n-\beta}}\,\frac{d\varrho}{\varrho}\,,\qquad \qquad \beta \in (0,n]
$$
with $B(x_0,\varrho)$ being the open ball centered at $x_0$, with radius $\varrho$.
In fact, it is not difficult to see that
\eqn{equivalenza}
$$
\pot{\beta}(x_0,R)\lesssim \int_{B_R(x_0)}\frac{d \mu(y)}{|x_0-y|^{n-\beta}}=I_\beta(\mu\llcorner B(x_0,R))(x_0)\leq I_\beta(\mu)(x_0)
$$
holds provided $\mu$ is a non-negative measure. A question is now, {\em is it possible to give an
analogue
of estimates \trif{stima0} in the case of general quasilinear equations such as for instance, the degenerate $p$-Laplacean equation}
\eqn{plap}
$$
-\divo\, (|Du|^{p-2}Du)=\mu \, ?
$$
A first answer has been given in the papers \cite{KM, TW}, where - for suitably defined solutions to \rif{plap} - the authors prove the following pointwise {\em zero order estimate} - i.e. for $u$ - when $p \leq n$,
via non-linear Wolff potentials:
\eqn{KMm}
$$
|u(x_0)| \leq c \left( \mean{B(x_0,R)}|u|^{p-1}\, dx\right)^{\frac{1}{p-1}}  + c\ww(x_0,2R)\,,
$$
where the constant $c$ depends on the quantities $n,p$, and
\eqn{newpot}
$$
{\bf W}^{\mu}_{\beta, p}(x_0,R):= \int_0^R \left(\frac{|\mu|(B(x_0,\varrho))}{\varrho^{n-\beta p}}\right)^{\frac{1}{p-1}}\, \frac{d\varrho}{\varrho}\qquad \qquad  \beta \in (0,n/p]\,,
$$
is the non-linear Wolff potential of $\mu$. Of course we are here using the standard notation concerning integral averages
$$
\mean{B(x_0,R)}|u|^{q}\, dx := \frac{1}{|B(x_0,R)|} \int_{B(x_0,R)}|u|^{q}\, dx\,.
$$
Estimate \rif{KMm}, which extends to a whole family of general quasi-linear equations, and which is commonly considered as a basic result in the theory of quasi-linear equations, is the natural non-linear analogue
of the first linear estimate appearing in \rif{stima0}. {\em Here we present the non-linear analogue
of the second estimate in \trif{stima0}, thereby giving a pointwise gradient estimate via non-linear potentials which upgrades \trif{newpot} up to the gradient/maximal level.}

\section{Degenerate Elliptic estimates}
In this section the growth exponent $p$ will be a number such that $p\geq 2$, we shall therefore treat possibly degenerate elliptic equations when $p\not=2$.
Specifically, we shall consider general non-linear, possibly degenerate equations with $p$-growth of the type
\eqn{baseq}
$$
-\divo \ a(x,Du)=\mu\,.
$$
whenever $\mu$ is a Radon measure with finite total mass defined on $\Omega$; eventually letting $\mu(\er^n \setminus \Omega)=0$, without loss of generality we may assume that $\mu$ is defined on the whole $\er^n$. 
The continuous vector field $a \colon \Omega \times \er^n \to \er^n$ is assumed to be $C^1$-regular in the gradient variable $z$, with $a_z(\cdot)$ being Carath\'eodory regular and
satisfying the following {\em growth, ellipticity and continuity assumptions}:
\eqn{asp}
$$
\left\{
    \begin{array}{c}
    |a(x,z)|+|a_{z}(x,z)|(|z|^2+s^2)^{\frac{1}{2}} \leq L(|z|^2+s^2)^{\frac{p-1}{2}} \\ [3 pt]
    \nu^{-1}(|z|^2+s^2)^{\frac{p-2}{2}}|\lambda|^{2} \leq \langle a_{z}(x,z)\lambda, \lambda
    \rangle\\[3pt]
    |a(x,z)-a(x_0,z)|\leq L_1\omega(|x-x_0|)(|z|^2+s^2)^{\frac{p-1}{2}}\,,
    \end{array}
    \right.
$$
whenever $x,x_0 \in \Omega$ and $z, \lambda \in \er^n$, where $0< \nu\leq  1 \leq L$ and $s\geq 0, L_1\geq 1$ are fixed parameters. When $p>2$ we shall assume that there exists a positive $\alpha< \min\{1, p-2\}$ such that the H\"older continuity property
\eqn{asp22}
$$
|a_{z}(x,z_2)-a_{z}(x,z_1)|\leq L|z_2-z_1|^{\alpha}(|z_1|^2+|z_2|^2+s^2)^{\frac{p-2-\alpha}{2}}
$$
holds whenever $z_1, z_2 \in \er^n$ and $x \in \Omega$. Here $\omega\colon [0,\infty) \to [0,\infty)$ is a modulus of continuity i.e. a non-decreasing function such that $\omega(0)=0$ and $\omega(\cdot)\leq 1$.
On such a function we impose a natural decay property, which is essentially optimal for the result we are going to have, and prescribes a {\em Dini continuous dependence of the partial map} $x \mapsto a(x,z)/(|z|+s)^{p-1}$:
\eqn{intdini}
$$
\int_0^R [\omega(\varrho)]^{\frac{2}{p}}\,\frac{ d \varrho}{\varrho} := d(R)< \infty\,,
$$
for some $R>0$. The prototype of \rif{baseq} is - choosing $s=0$ and omitting the $x$-dependence - clearly given by the $p$-Laplacean equation
\rif{plap}, which satisfies \rif{asp22} whenever $\alpha< \min\{1, p-2\}$.
In the following, when a measure $\mu$ actually turns out to be an $L^1$-function, we shall use the standard notation
$$
|\mu|(A):=\int_{A} |\mu(x)|\, dx\,,
$$
whenever $A$ is a measurable set on which $\mu$ is defined.

 In this paper we shall present our results in the form of a priori estimates - i.e. when solutions and data are taken to be more regular than needed, for instance $u \in C^{1}(\Omega)$  and $\mu \in L^1(\Omega)$  - but they actually hold, via a standard approximation argument, for general weak and very weak solutions - i.e. distributional solutions which are not in the natural space $W^{1,p}(\Omega)$ - to measure data problems such as, for instance
\eqn{Dir1}
$$
\left\{
    \begin{array}{cc}
    -\divo \ a(x,Du)=\mu & \qquad \mbox{in $\Omega$}\\
        u= 0&\qquad \mbox{on $\partial\Omega$\,,}
\end{array}\right.
$$
where $\mu$ is a general Radon measure with finite total mass, defined on $\Omega$. The reason for such a
choice is that the approximation argument in question leads to different notions of solutions, according to the regularity/integrability properties of the right hand side $\mu$. We do not want to enter in such details too much, for which we refer to \cite{DDM, mis3}, and therefore we confine ourselves to the neat a priori estimate form of the results.

For instance, in the case \rif{Dir1} with $\mu$ being a genuinely Radon measure, in \cite{DDM, mis3} we consider the so called Solutions Obtained by Limit of Approximations (SOLA), which is a standard class considered when dealing with measure data problems. Such solutions are in particular unique in the case $p=2$, as proved in \cite{boccardo, TW2}. Finally, if the right hand side of \rif{baseq} is integrable enough to deduce that $\mu \in W^{-1,p'}(\Omega)$, then our results apply to general weak energy solutions $u \in W^{1,p}(\Omega)$ to \rif{baseq}.

The first result we present is  now
\begin{theorem}[Non-linear potential gradient bound]\label{mainx} Let $u \in C^{1}(\Omega)$, be a weak solution to \trif{baseq} with $\mu \in L^1(\Omega)$, under the assumptions \trif{asp}. Then there exists a constant $c \equiv c (n,p,\ratio, L_1, \alpha)>1$, and a positive radius $R_0$ depending only on $n,p,\ratio,L_1, \omega(\cdot), \alpha$, such that
the pointwise estimate
\eqn{mainestx}
$$
|Du(x_0)| \leq c\left(\mean{B(x_0,R)}(|Du|+s)^{\frac{p}{2}}\, dx\right)^{\frac{2}{p}}  + c\npma(x_0,2R)
$$
holds whenever $B(x_0,2R)\subseteq \Omega$, and $R\leq R_0$. Moreover, when the vector field $a(\cdot)$ is independent of $x$ - and in particular for the $p$-Laplacean operator \trif{plap} - estimate \trif{mainestx} holds with no restriction on $R$.
\end{theorem}
The potential $\npma$ appearing in \rif{mainestx} is the natural one since its shape respects the scaling properties of the equation with respect to the estimate in question; compare with the linear estimates \rif{stima0}. When extended to general weak solutions {\em estimate \trif{mainestx} tells us the remarkable fact that the boundedness of $Du$ at a point $x_0$ is independent of the solution $u$, and of the vector field $a(\cdot)$ considered, but only depends on the behavior of $|\mu|$ in a neighborhood of $x_0$.}

A particularly interesting situation occurs in the case $p=2$, when we have a pointwise potential estimate which is completely similar to the second one in \rif{stima0}, and that we think deserves a statement of its own, that is
\begin{theorem}[Linear potential gradient bound]\label{mainx2} Let $u \in C^{1}(\Omega)$, be a weak solution to \trif{baseq} with $\mu \in L^1(\Omega)$,
under the assumptions \trif{asp} considered with $p=2$. Then there exists a constant $c \equiv c (n,p,\ratio, L_1)>0$, and a positive radius $R_0 \equiv R_0(n,p,\ratio,L_1, \omega(\cdot))$ such that
the pointwise estimate
\eqn{mainest}
$$
|D u(x_0)| \leq c \mean{B(x_0,R)}(|D u|+s)\, dx  + c\potm{1}(x_0,2R)
$$
holds whenever $B(x_0,2R)\subseteq \Omega$, and $R\leq R_0$.
Moreover, when the vector field $a(\cdot)$ is independent of the variable $x$, estimate \trif{mainest} holds with no restriction on $R$.
\end{theorem}
Beside their intrinsic theoretical interest, the point in estimates \rif{mainestx}-\rif{mainest} is that they allow to unify and recast essentially all the gradient $L^q$-estimates for quasilinear equations in divergence form; moreover they allow for an immediate derivation of estimates in intermediate spaces such as interpolation spaces. We refer to the recent survey \cite{mm} for an account of such estimates. Indeed, by \rif{mainestx} it is clear that the behavior of $Du$ can be controlled by that $\npma$, which is in turn known via the behavior of Riesz potentials. In fact, this is a consequence of the pointwise bound of the Wolff potential via the Havin-Maz'ja non linear potential \cite{AdMe, MH, AdHe}, that is
\eqn{rieszbound}
$$
\npma(\cdot,\infty) = \int_{0}^\infty \left(\frac{|\mu|(B(x_0,\varrho))}{\varrho^{n-1}}\right)^{\frac{1}{p-1}} \frac{d \varrho}{\varrho} \leq c I_{\frac{1}{p}}\left\{\left[I_{\frac{1}{p}}(|\mu|)\right]^{\frac{1}{p-1}}\right\}(x_0)\,.
$$
Ultimately, thanks to \rif{rieszbound} and to the well-known properties of the Riesz potentials, we have
\eqn{wolff}
$$
 \mu \in L^q  \Longrightarrow \npma \in L^{\frac{nq(p-1)}{n-q}} \qquad  q \in (1,n)\,,
$$
while Marcikiewicz spaces must be introduced for the borderline case $q=1$. Inequality \rif{wolff} immediately allows to recast the classical gradient estimates for solutions to \rif{Dir1} such as those due to Boccardo \& Gall\"ouet \cite{BG1, BG2} - when $q$ is \ap small" - and Iwaniec \cite{I} and DiBenedetto \& Manfredi \cite{DM} - when $q$ is \ap large" -  that is, for solutions to \rif{Dir1} it holds that
$$
 \mu \in L^q  \Longrightarrow Du \in L^{\frac{nq(p-1)}{n-q}} \qquad  q \in (1,n)\,.
$$
Moreover, since the operator $\mu \mapsto \npma$ is obviously sub-linear, using the estimates related to \rif{wolff} and classical interpolation theorems for sub-linear operators one immediately gets estimates in refined scales of spaces such Lorentz or Orlicz spaces, recovering some estimates of Talenti \cite{Talenti}, but directly for the gradient of solutions, rather than for solutions themselves.

Another point of Theorem \ref{mainx} is that it allows to prove an essentially optimal Lipschitz continuity criterium with respect to the regularity of coefficients \rif{intdini}, that is
\eqn{minass}
$$\npma(\cdot,R) \in L^{\infty}(\Omega), \ \mbox{for some}\  R>0 \Longrightarrow Du \in L^{\infty}_{\loc}(\Omega, \er^n)\,,$$
and moreover the local bound
\eqn{minimalcri}
$$
\|Du\|_{L^{\infty}(B_{R/2})} \leq c \left(\mean{B(x_0,R)}(|Du|+s)^{\frac{p}{2}}\, dx\right)^{\frac{2}{p}} + c\left\|\npma(\cdot,R)\right\|_{L^{\infty}(B_{R})}
$$
holds whenever $B_{2R}\subseteq \Omega$.

We finally recall that another consequence of the classical estimate \rif{rieszbound} and of \rif{mainestx} is
\eqn{mainestxr}
$$
|Du(x_0)| \leq c\left(\mean{B(x_0,R)}(|Du|+s)^{\frac{p}{2}}\, dx\right)^{\frac{2}{p}}  + cI_{\frac{1}{p}}\left\{\left[I_{\frac{1}{p}}(|\mu|)\right]^{\frac{1}{p-1}}\right\}(x_0)\,,
$$
which holds whenever $B(x_0,2R)\subset \Omega$ satisfies the conditions imposed in Theorem \ref{mainx}. Here we recall the reader that we have previously extended $\mu$ to the whole space $\er^n$.
\section{Parabolic first, and zero order estimates} Our aim here is not only to give a parabolic version of the elliptic estimate \rif{mainestx}, but also to give a zero order estimate, that is the parabolic analog of the zero order elliptic estimate \cite{KM}, the validity of which {\em was yet considered to be an open issue.} We consider quasilinear parabolic equations of the type
\eqn{basicpar}
$$
u_t - \divo \ a(x,t,Du)= \mu\,,
$$
in a cylindrical domain $\Omega_T:= \Omega\times (-T,0)$,
where as in the previous section $\Omega \subset \er^n$, $n \geq 2$ and $T>0$.
The vector-field $a\colon\Omega_T\times\er^n\to\er^n$ is assumed to be
Carath\`eodory regular together with $a_z(\cdot)$, and indeed being $C^1$-regular with respect to the gradient variable $z \in \er^n$, and satisfying the following standard growth, ellipticity/parabolicity and continuity conditions:
\begin{equation}\label{par1}
    \left\{
    \begin{array}{c}
 |a(x,t,z)|+|a_{z}(x,t,z)|(|z|+s) \leq L(|z|+s) \\[4pt]
    \nu|\lambda|^2 \leq \langle a_{z}(x,t,z)\lambda, \lambda\rangle  \\[4pt]
    |a(x,t,z)-a(x_0,t,z)|\leq L_1\omega(|x-x_0|)(|z|+s)
    \end{array}
    \right.
\end{equation}
for every choice of $x, x_0\in\Omega$, $z,\lambda \in\er^n$ and $ t \in (-T,0)$; here the function $\omega\colon \er^+ \to \er^+$ is as in \rif{asp}$_3$. Note that anyway we are assuming no continuity on the map $t \mapsto a(\cdot, t, \cdot)$, which is considered to be a priori only measurable. In other words we are considering the analog of assumptions \rif{asp} for $p=2$; the reason we are adopting this restriction is that when dealing with the evolutionary $p$-Laplacean operator estimates assume the usual form only when using so called \ap intrinsic cylinders", according the parabolic $p$-Laplacean theory developed by DiBenedetto \cite{D}. These are - unless $p=2$ when they reduce to the standard parabolic ones - cylinders whose size locally depends on the size of the solutions itself, therefore a formulation of the estimates via non-linear potentials - whose definition is built essentially using a standard family of balls and it is therefore \ap universal" - is not immediate and will be the object of future investigation. We refer to \cite{AM} for global gradient estimates. 

In order to state our results we need some additional terminology. Let us recall that given points $(x,t), (x_0,t_0)\in\er^{n+1}$ the standard parabolic metric is defined by
\eqn{parme}
$$
d_{\rm par}((x,t),(x_0,t_0)):=\max\{|x-x_0|, \sqrt{|t-t_0|}\} \thickapprox \sqrt{|x-x_0|^2+|t-t_0|}
$$
and the related metric balls with radius $R$ with respect to this metric are given by cylinders $B(x_0,R)\times (t_0-R^2, t_0+R^2)$. The \ap caloric" Riesz potential - compare with elliptic one defined in \rif{ellr}, and with \cite{AB}, for instance - is now built starting from \rif{parme}
\begin{equation}\label{parr}
    I_\beta (\mu)((x,t)):=\int_{\er^{n+1}}\frac{d\mu((\tilde x,\tilde t ))}{d_{\rm par}((\tilde x,\tilde t ),(x,t))^{N-\beta}}\, ,
    \qquad 0 < \beta \leq N:=n+2\,,
\end{equation}
whenever $(x,t) \in \er^{n+1}$. In order to be used in estimates for parabolic equations, it is convenient to introduce its local version
via the usual backward parabolic cylinders - with \ap vertex" at $(x_0,t_0)$ - that is
\eqn{caloricc}
$$
Q(x_0,t_0;R):= B(x_0,R)\times (t_0-R^2, t_0)\,,
$$
so that we define
\eqn{calpot}
$$
 {\bf I}_{\beta}^\mu(x_0,t_0;R):=
 \int_0^R \frac{\mu(Q(x_0,t_0;\varrho))}{\varrho^{N-\beta}}\, \frac{d\varrho}{\varrho}\qquad \mbox{where}\ \beta \in (0,N]\,.
$$
The main result in the parabolic case is
\begin{theorem}[Parabolic potential gradient bound]\label{mainp} Under the assumptions \trif{par1} and \trif{intdini}, let $u\in C^0(-T,0; L^2(\Omega))$ be a weak solution to \trif{basicpar} with $\mu \in L^\infty(\Omega_T)$ and such that $Du\in C^0(\Omega_T)$. Then there exists a constant $c\equiv c(n,\nu, L) $ and a radius $R_0\equiv R_0(n,\nu, L, L_1, \omega(\cdot)) $ such that the following estimate:
\eqn{parest1}
$$
|D u(x_0,t_0)|  \leq c \mean{Q(x_0,t_0;R)}(|D u|+s)\, dx\, dt + c{\bf I}_{1}^{|\mu|}(x_0,t_0;2R) \,,
$$
holds whenever $Q(x_0,t_0; 2R)\subseteq \Omega$, and $R \leq R_0$. When the vector field $a(\cdot)$ is independent of the space variable $x$, estimate \trif{parest1} holds with no restriction on $R$.
\end{theorem}
Again, as in the elliptic case, estimate \rif{parest1} also holds for solutions to general measure data problems as
\begin{equation}\label{parpro}
    \left\{
    \begin{array}{cc}
    u_t-\divo \, a(x,t,Du)=\mu&\mbox{in $ \Omega_T$}\\[3pt]
    u=0&\mbox{on $\partial_{\rm par}\Omega_T$\,,}
    \end{array}
    \right.
\end{equation}
where $\mu$ is a general Radon measure with finite mass on $\Omega_T$, that we shall
again consider to be defined in the whole $\er^{n+1}$. In the spirit of the elliptic result \rif{minimalcri} we have the following implication, which provides a boundedness criteria for the spatial gradient, under the Dini continuity assumption for the spatial coefficients
stated in \rif{intdini}:
\eqn{minasspar}
$${\bf I}_{1}^{|\mu|}(\cdot;R) \in L^{\infty}(\Omega_T), \ \mbox{for some}\  R>0 \Longrightarrow Du \in L^{\infty}_{\loc}(\Omega_T, \er^n)\,.$$
We conclude with the zero order potential estimate, which applies to general equations of the type \rif{basicpar} when considered with a measurable dependence upon the coefficients $(x,t)$. The relevant hypotheses here are the following standard growth and monotonicity properties:
\begin{equation}\label{par2}
    \left\{
    \begin{array}{c}
 |a(x,t,z)|\leq L(|z|+s) \\[4pt]
    \nu|z_2-z_1|^2 \leq \langle a(x,t,z_2)-a(x,t,z_1), z_2-z_1\rangle
    \end{array}
    \right.
\end{equation}
which are assumed to hold whenever $(x,t)\in \Omega_T$ and $z,z_1,z_2 \in \er^n$. In particular, since the pointwise bound will be derived on $u$, rather than on $Du$, we do not need any differentiability assumption on $a(\cdot)$ with respect to the spatial gradient variable $z$-variable, assumptions \rif{par2} are clearly weaker than \rif{par1}.
\begin{theorem}\label{mainp2} Under the assumptions \trif{par2}, let $u \in  L^2(-T, 0;W^{1,2}(\Omega))\cap C^0(\Omega_T)$ be a weak solution to \trif{basicpar} with $\mu \in L^1(\Omega_T)$. Then there exists a constant $c$, depending only on $n,\ratio, L_1$ such that the following inequality holds whenever $Q(x_0,t_0; 2R)\subseteq \Omega$:
\eqn{lastpest}
$$
|u(x_0,t_0)| \leq c \mean{Q(x_0,t_0;R)}(|u|+s)\, dx\, dt  +  c{\bf I}_{2}^{|\mu|}(x_0,t_0;2R) + cRs \,.
$$
\end{theorem}
\section{A non-local Caccioppoli's inequality}
In \cite{DDM, mis3} we have developed more than one approach to the proof of the pointwise gradient estimates
via non-linear potentials. Here we shall present one of these,  taken form \cite{mis3}, for the case $p=2$, and for simplicity restricting to equations 
with no coefficients i.e. of the type
\eqn{basis}
$$
\divo\, a(Du)=\mu\,.
$$
{\em We believe that such method of proof is of independent technical interest} since it potentially applies to all those problems with a lack of full differentiability, as it will be clear in a few lines. Moreover, we shall see that in the case \rif{basis} estimate \rif{mainest} holds component-wise; see
\rif{lastest} below. The assumptions considered for \rif{basis} are of course
\eqn{aspdue}
$$
    \nu|\lambda|^{2} \leq \langle a_z(z)\lambda, \lambda
    \rangle\,,\qquad \qquad | a_z(z)| \leq L\,,\qquad \qquad |a(0)|\leq L\,.
$$
which hold whenever $z, \lambda \in \er^n$, where $0 < \nu \leq L$. The presentation of this
technique is indeed one of the objectives of \cite{mis3}. Aiming at the explanation of a general viewpoint, let us recall that for energy solutions $u \in W^{1,2}(\Omega)$ to homogeneous equations of the type
\eqn{hom}
$$
\divo\, a(Du)=0
$$
the local boundedness of the gradient is achieved by {\em first} differentiating the equation \rif{hom}, proving that $Du \in W^{1,2}_{\loc}(\Omega)$, and {\em then} observing that $v:= D_\xi u$ solves the linear equation with measurable coefficients
$$
\divo (A(x)Dv)=0 \qquad \qquad A(x):= a_z(Du(x))\,.
$$
At this stage the boundedness of $D_\xi v$ follows applying an iteration method, as for instance the one devised in the pioneering work of DeGiorgi \cite{DG}. This is
in turn
based on the use of {\em Caccioppoli's inequalities on level sets}, that is, denoting
$$
(w-k)_+:= \max \{w-k,0\}\,,\qquad  \qquad (w-k)_-:= \max \{k-w,0\} $$
we have that inequalities of the type
\eqn{cacc1}
$$
\int_{B_{R/2}} |D(D_\xi u-k)_+|^2\, dx  \leq \frac{c}{R^2}\int_{B_{R/2}} |(D_\xi u-k)_+|^2\, dx
$$
and similar variants, for instance involving $(D_\xi u-k)_-$, hold whenever $k \in \er$. In turn, the iteration of such inequalities yields the boundedness of $D_\xi u$.
In such an iteration, {\em one controls the level sets of $D_{\xi} u$} via the higher order derivatives $D(D_\xi u-k)_+$ and Sobolev embedding theorem, building a geometric iteration in which, at every step, the gain is dictated by the Sobolev embedding exponent.

Applying such a reasoning to the case \rif{basis} seems to be difficult, as even in the simplest case \rif{poisson} it is in general false that $Du \in W^{1,1}(\Omega)$ when the right hand side $\mu$ is just a measure, or an $L^1$-function. On the other hand, a result of \cite{mis1} states that, although Calder\'on-Zygmund theory does not apply in the classical $W^{1,1}$-sense, when considering the borderline case when $\mu$ is a measure or lies
in $L^1$, it nevertheless holds {\em provided the right functional setting} is considered, i.e. using Fractional Sobolev spaces.
Indeed, for SOLA to measure data problems as \rif{basis} it holds that
\eqn{fCZ}
$$
Du  \in W^{1-\ep, 1}_{\loc}(\Omega, \er^n)\qquad \qquad \mbox{for every}\ \ep \in (0,1)\,,
$$
with related explicit a priori local estimates; see \cite[Theorem 1.2]{mis1} for precise statements.
We here recall that, for a bounded open set $A
\subset \er^n$ and $k \in \en$, parameters $\alpha \in (0,1)$ and $q
\in [1,\infty )$, the fractional Sobolev space $W^{\alpha ,q}(A,\er^k )$ consists of those measurable mappings $w\colon \Omega \to \er^k$ such that
the following Gagliardo-type norm is finite:
\begin{eqnarray}
\nonumber \| w \|_{W^{\alpha,q}(A )} & := & \left(\int_{A}
|w(x)|^{q}\, dx \right)^{\frac{1}{q}} + \left(\int_{A} \int_{A}
\frac{|w(x)
- w(y) |^{q}}{|x-y|^{n+\alpha q}} \ dx \, dy \right)^{\frac{1}{q}}\\
&   =:&  \|w \|_{L^q(A)} +[ w ]_{\alpha ,q;A} <
\infty\;.\label{defrazionari}
\end{eqnarray}
With such a notation \rif{fCZ} means that
\eqn{fCZ2}
$$
[ Du ]_{1-\ep ,1;\Omega'}= \int_{\Omega'} \int_{\Omega'}
\frac{|Du(x)
- Du(y) |}{|x-y|^{n+1-\ep}} \ dx \, dy < \infty
$$
holds for every $\ep \in (0,1)$, and every subdomain $\Omega' \Subset \Omega$; the previous quantity is intuitively the $L^1$-norm of the \ap $(1-\ep)$-order derivative" of $Du$, roughly denotable by $D^{1-\ep}Du$. The inequality in \rif{fCZ2} let us think that Caccioppoli type inequality \rif{cacc1} should be replaced by a fractional order
version, and using the $L^1$-norm, rather than the $L^2$-one. Indeed we have the following theorem, that we again for simplicity state under the form of a priori estimate - i.e. assuming more regularity $u\in W^{1,2}(\Omega)$ and $\mu \in L^2(\Omega)$ (this can be again removed via an approximation scheme, and by considering suitable definitions of solutions). Needless to say, what it matters here is the precise form of the a priori estimate.
 \begin{theorem}[Non-local Caccioppoli inequality]\label{diffra} Let $u\in W^{1,2}(\Omega)$ be a weak solution to \trif{basis} with $\mu \in L^2(\Omega)$, under the assumptions \trif{aspdue}; whenever $\xi \in \{1,\ldots,n\}$, $k \geq 0$, and whenever $B_{R} \subseteq \Omega$ is a ball with radius $R$, the inequality
\eqn{sob01}
$$
 [(|D_\xi u|-k)_+]_{\sigma, 1;B_{R/2}}  \leq  \frac{c}{R^{\sigma}} \int_{B_R} (|D_{\xi} u|-k)_+ \, dx + \frac{cR|\mu|(B_R)}{R^{\sigma}}\,,
$$
holds for every $\sigma< 1/2$, where the constant $c$ depends only on $n,\ratio,\sigma$.
\end{theorem}
Comparing \rif{sob01} and \rif {cacc1}, {\em Theorem \ref{diffra} tells us that for quasilinear equations Caccioppoli's inequalities are a robust tool that keeps holding at intermediate derivatives/integrability levels.}
We do think that the idea of using non-local Caccioppoli inequalities instead of the usual ones is interesting in itself as it leads to certain types of iterations which work without fully differentiating the equation; in turn, this could apply to all those problems with a lack of full differentiability. We indeed explicitly note here that a fractional Caccioppoli inequality has been indeed derived for notwithstanding the problems has integer order. The proof of the
inequality is developed in \cite{mis3} and has as a starting point some techniques introduced in \cite{KMin, mis1}.

The idea is now rather natural: inequality \rif{sob01} serves to start an iteration in which, at each stage we control the level set of $D_{\xi} u$ via the fractional derivative $D^\sigma (D_\xi u)$ and the fractional version of Sobolev embedding theorem. We come up again with a geometric iteration whose step is in turn dictated by the fractional Sobolev embedding exponent. A point we want to emphasize, is that, as clearly inferrable from \cite{mis3}, inequality \rif{sob01} contains all the information about the pointwise gradient estimate, no matter how small $\sigma$ is taken. As a matter of fact in the following we are not using explicitly the fact that $u$ is a solution, but rather the fact that $D_\xi u$ satisfies \rif{sob01}. For this reason, we shall report the next result in an abstract way. Moreover, we think that the formulation below could be useful in different contexts.
 \begin{theorem}[De Giorgi's fractional iteration]\label{degfra} Let $w\in L^1(\Omega)$ be a function with the property that there exist $\sigma \in (0,1)$ and $c_1\geq 1$, and a Radon measure $\mu$, such that whenever $B_{R} \subseteq \Omega$ is a ball with radius $R$ and $k \geq 0$, the inequality
\eqn{sob22}
$$
 [(|w|-k)_+]_{\sigma, 1;B_{R/2}}  \leq  \frac{c_1}{R^{\sigma}} \int_{B_R} (|w|-k)_+ \, dx + \frac{c_1R|\mu|(B_R)}{R^{\sigma}}\,,
$$
holds. Then the following estimate:
\eqn{sob33}
$$
|w (x_0)| \leq c \mean{B(x_0,R)}|w|\, dx  + c\potm{1}(x_0,2R)
$$
holds whenever $B(x_0,2R)\subset \Omega$, where the constant $c$ depends on $c_1, n, \sigma$.
\end{theorem}
The dependence of the constant $c$ appearing in \rif{sob33} is not surprisingly as follows:
$$
\lim_{\sigma \to 0} c = \infty \qquad \mbox{and}\qquad \lim_{c_1 \to \infty} c = \infty\,.
$$
Now we just have to conclude merging the last two theorems. Indeed, by Theorem \ref{diffra} we have assumption \rif{sob22} from Theorem \ref{degfra} satisfied by $w \equiv D_\xi u$. In turn, applying Theorem \ref{degfra} with such a choice of $w$ we conclude with the desired pointwise gradient bound
\eqn{lastest}
$$
|D_\xi u(x_0)| \leq c \mean{B(x_0,R)}|D_\xi u|\, dx  + c\potm{1}(x_0,2R)\,.
$$
The last estimate clearly implies \rif{mainest}, being actually stronger since it holds for each single component of the gradient.

{\bf Acknowledgments.} This research is supported by the ERC grant 207573 ``Vectorial Problems", and by MIUR via the national project \ap Calcolo delle Variazioni".


\begin{thebibliography}{99}


 \bibitem {AM}\name[Acerbi, E.]\et
\name[Mingione, G.]: Gradient estimates for a class of parabolic
systems. {\em Duke Math. J.} 136 (2007), 285--320.


\bibitem {AB} \name[Adams, D.~R.]\et \name[Bagby, R.~J.]: Translation-dilation invariant estimates for Riesz potentials. {\em Indiana Univ.~Math.~J.} 23 (1974), 1051--1067.

\bibitem {AdHe} \name[Adams, D.~R.]\et \name[Hedberg, L.~I.]: {\em Function spaces and potential
theory.}  Grundlehren der Mathematischen Wissenschaften
 314.
Springer-Verlag, Berlin, 1996.

\bibitem {AdMe} \name[Adams, D.~R.]\et \name[Meyers, N.~G.]: Thinnes and Wiener criteria for non-linear potentials. {\em Indiana Univ.~Math.~J.} 22 (1972), 169--197.


\bibitem {Ad}\name[Adams, R.A.]: {\em  Sobolev Spaces.} Academic
Press, New York, 1975.



\bibitem {boccardo} \name[Boccardo, L.]: Problemi differenziali ellittici e parabolici con dati misure [Elliptic and parabolic differential problems
with measure data].
{\em Boll.~Un.~Mat.~Ital.~A (7)} 11 (1997),
 439--461.


\bibitem {BG1}
\name[Boccardo, L.]\et \name[Gallou\"et, T.]:
 Nonlinear elliptic
and parabolic equations involving measure data. {\em
J.~Funct.~Anal.~}87 (1989), 149--169.


\bibitem {BG2}
\name[Boccardo, L.]\et \name[Gallou\"et, T.]: Nonlinear elliptic
equations with right-hand side measures. {\em Comm.~Partial
Differential Equations} 17 (1992), 641--655.

\bibitem {D} \name[DiBenedetto, E.]: {\em Degenerate parabolic equations.} Universitext. Springer-Verlag, New
York, 1993.

\bibitem {DM} \name[DiBenedetto, E.]\et\name[Manfredi, J.J.]: On the higher
integrability of the gradient of weak solutions of certain
degenerate elliptic systems. {\em Amer.~J.~Math.} 115 (1993),
1107--1134.




\bibitem {DG} \name[De Giorgi, E.]: Sulla differenziabilit\`a e l'analiticit\`a delle
 estremali degli integrali multipli regolari. {\em
 Mem.~Accad.~Sci.~Torino Cl.~Sci.~Fis.~Mat.~Nat.~(III)} 125 3 (1957), 25--43.


 \bibitem {DDM} \name[Duzaar, F.]\et\name[Mingione, G.]:
Gradient estimates via non-linear potentials. {\em Submitted 2009}.








\bibitem {GW} \name[Gr\"uter, M.]\et\name[Widman, K.O.]: The Green function for uniformly elliptic equations. {\em manuscripta math.} 37 (1982), 303--342.


    \bibitem {MH} \name[Havin, M.]\et\name[Maz'ja, V. G.]:
A nonlinear potential theory. {\em Russ.~Math.~Surveys} 27 (1972), 71--148.

\bibitem {HW} \name[Hedberg, L.I.]\et \name[Wolff, T.]:
Thin sets in nonlinear potential theory.
{\em Ann.~Inst.~Fourier (Grenoble)} 33 (1983), 161--187.


\bibitem {I} \name[Iwaniec, T.]: Projections onto gradient fields and
$L\sp{p}$-estimates for degenerated elliptic operators. {\em Studia
Math.} 75 (1983), 293--312.



\bibitem {KM} \name[Kilpel\"ainen, T.]\et\name[Mal\'y, J.]: The Wiener test and potential
estimates for quasilinear elliptic equations. {\em Acta Math.}~172
(1994), 137--161.


\bibitem {KMin} \name[Kristensen, J.]\et\name[Mingione, G.]:
The singular set of minima of integral functionals. {\em
Arch.~Ration.~Mech.~Anal.} 180 (2006), 331--398.




\bibitem {mis1} \name[Mingione, G.]: The Calder\'on-Zygmund theory for elliptic problems with measure data.
{\em Ann~Scu.~Norm.~Sup.~Pisa Cl.~Sci.~(5)} 6 (2007),
195--261.


\bibitem {mis2} \name[Mingione, G.]: Gradient estimates below the duality exponent.
{\em Submitted} (2008).

\bibitem {mm} \name[Mingione, G.]: Towards a non-linear Calder\'on-Zygmund theory. {\em Quaderni di Matematica (Alvino, Mercaldo, Murat \&  Peral eds.)} (2009).


\bibitem {mis3} \name[Mingione, G.]: Gradient potential estimates.
{\em Submitted} (2008).




\bibitem{Talenti} \name[Talenti, G.]: Elliptic equations and rearrangements. {\em Ann~Scu.~Norm.~Sup.~Pisa Cl.~Sci.~(4)} 3 (1976), 697--717.

\bibitem{TW} \name[Trudinger, N.S.]\et\name[Wang, X.J.]: On the weak continuity of
elliptic operators and applications to potential theory. {\em
Amer.~J.~Math.} 124 (2002), 369--410.



\bibitem{TW2} \name[Trudinger, N.S.]\et\name[Wang, X.J.]: Quasilinear elliptic equations with signed measure data.
{\em
Disc.~Cont.~Dyn.~Systems A} 23 (2009), 477--494.



\end{thebibliography}
\end{document}